\documentclass[11pt,a4paper]{amsart}
\numberwithin{equation}{section}
\usepackage{rotating}
\usepackage{mathrsfs}
\usepackage{amscd,amssymb,amsopn,amsmath,amsthm,graphics,amsfonts,enumerate,verbatim,calc}
\usepackage[all]{xy}
\usepackage[utf8]{inputenc}
\usepackage{color}

\usepackage[a4paper,margin=1in]{geometry}
\usepackage[utf8]{inputenc}
\usepackage{amscd,amssymb,amsopn,amsmath,amsthm,graphics,amsfonts,enumerate,verbatim,calc}
\usepackage[T1]{fontenc}
\usepackage{lmodern}
\usepackage{enumitem}
\usepackage{hyperref}
\usepackage{microtype}\usepackage{hyperref}
\usepackage{geometry}
\usepackage{enumitem}
\usepackage{booktabs}
\usepackage{tikz-cd}
\hypersetup{colorlinks=true,linkcolor=blue,citecolor=blue,urlcolor=blue}
\newtheorem{theorem}{Theorem}[section]
\newtheorem{proposition}[theorem]{Proposition}
\newtheorem{corollary}[theorem]{Corollary}
\newtheorem{lemma}[theorem]{Lemma}
\newtheorem{fact}[theorem]{Fact}
\newtheorem{example}[theorem]{Example}\newtheorem{notation}[theorem]{Notation}

\newtheorem{remark}[theorem]{Remark}
\newtheorem{question}[theorem]{Question}

\newtheorem{definition}[theorem]{Definition}

\newcommand{\ulim}{\operatorname*{ulim}}
\newcommand{\U}{\mathcal{U}}

\DeclareMathOperator{\edim}{edim}
\DeclareMathOperator{\length}{\ell}
\newcommand{\m}{\mathfrak m}
\newcommand{\n}{\mathfrak n}

\newcommand{\ol}[1]{\overline{#1}}

\newcommand{\gr}{\operatorname{gr}}

\newcommand{\ord}{\operatorname{ord}}

\newcommand{\e}{\mathrm e}

\newcommand{\mm}{\mathfrak m}
\newcommand{\ideal}{\mathfrak a}
\newcommand{\NN}{\mathbb Z_{\geq0}}
\newcommand{\one}{\mathbf 1}
\newcommand{\cl}[1]{\overline{#1}}
\DeclareMathOperator{\Int}{Int}
\DeclareMathOperator{\conv}{conv}

\begin{document}
	
\author[M. Asgharzadeh]{Mohsen Asgharzadeh}
\date{}

\address{M. Asgharzadeh}
\email{mohsenasgharzadeh@gmail.com}

	\title[Iterated socles and closure operations ]
{Iterated socles and closure operations }

\subjclass[2020]{ 13H10; 13B22; 13A50}
	\keywords{linkage; Goto number; ultraproducts; integral closure; socle; test ideals; mixed char. }

	\begin{abstract}
  Let $(R,\m)$ be a $d$-dimensional Cohen--Macaulay local ring, $d\geq2$ containing a field,
  and let $J\subseteq\m^s$ be a parameter ideal with $s\geq2$,
also, in the mixed case via some mild assumptions. We prove for that
  $J:\m^{(d-1)(s-1)}\subseteq\ol J$, and that
  $J:\m^{(d-1)(s-1)+1}\subseteq\ol J$ when $R$ is not regular.  We 
  show that Cohen--Macaulayness cannot be removed.  We also prove the characteristic-free exact
  formula $g(J)=(d-1)(s-1)+r$ for parameter ideals of order $s$ with a
  level tangent cone, giving singular mixed-characteristic Gorenstein and
  non-Gorenstein examples in which the singular bound is attained.
\end{abstract}

	\maketitle

	\section{Introduction}
	\label{sec:intro}

Let $(R,\mathfrak m)$ be a Noetherian local ring and let $J$ be an
$\mathfrak m$-primary ideal. The \emph{Goto number} of $J$ is
$ g(J):=\max\{t\geq 0\mid J:\mathfrak m^t\subseteq\overline J\}, $
where $\overline J$ denotes the integral closure of $J$.
Thus, the Goto number measures how far one can enlarge $J$ by
taking colons with powers of the maximal ideal while remaining
inside its integral closure.
Motivated by a result of Wang \cite{Wang}, Watanabe and Yoshida
established a remarkable lower bound for this invariant:
\begin{theorem}[Watanabe--Yoshida {\cite[Theorem 6.2]{WY}}]
	\label{thm:WY}
	Let $(R,\mathfrak m)$ be a $d$-dimensional Cohen--Macaulay local
	ring of characteristic $p>0$, where $d\geq 2$. Let $s\geq 2$
	and let $J\subseteq\mathfrak m^s$ be a parameter ideal. Then
	\begin{enumerate}[label=\textup{(\roman*)}]
		\item
		$(J:\mathfrak m^{(d-1)(s-1)})\subseteq\overline J;$
		\item if $R$ is not regular, then
		$	(J:\mathfrak m^{(d-1)(s-1)+1})\subseteq\overline J.$
	\end{enumerate}
\end{theorem}

In characteristic zero, Watanabe and Yoshida obtained related
results under additional hypotheses, including normality and
conditions on the canonical divisor. They also raised the question
of whether the characteristic assumption in the above theorem
could be removed. More precisely, they asked the following.

\begin{question}[See {\cite[Question 8.1]{WY}}]
	\label{q:charfree}
	Does the containment in Theorem~\ref{thm:WY} hold without any
	restriction on the characteristic?
\end{question}

They also asked whether the Cohen--Macaulay hypothesis could
be removed.

\begin{question}[See {\cite[Question 8.2]{WY}}]
	\label{q:nonCM}
	Can the Cohen--Macaulay hypothesis in Theorem~\ref{thm:WY}
	be removed?
\end{question}

A further question, posed in \cite[Question 1.2]{WY}, concerns
the determination of the Goto number itself. Related questions
about the containment of colon ideals in integral closures were
also considered in \cite{it}.
The purpose of this paper is to investigate these questions from
three complementary perspectives: characteristic-free
containments, mixed-characteristic methods, and exact
computations of Goto numbers.
Our first result is: 

\begin{theorem}
	\label{thm:field}
	Let $(R,\mathfrak m)$ be a Cohen--Macaulay local ring containing
	a field, with $d=\dim R\geq 2$. Let $s\geq 2$ and let
	$J\subseteq\mathfrak m^s$ be a parameter ideal. Then
	$(J:\mathfrak m^{(d-1)(s-1)})\subseteq\overline J.	$
	If $R$ is not regular, then also
	$(J:\mathfrak m^{(d-1)(s-1)+1})\subseteq\overline J.	$
\end{theorem}

The proof in characteristic zero uses the Lefschetz approximation
theory of Aschenbrenner and Schoutens. The basic idea is to
approximate a complete local ring of characteristic zero by
complete local rings of positive characteristic, apply the
Watanabe--Yoshida theorem to these approximating rings, and
descend the resulting containment by faithful flatness.

We next turn to mixed characteristic. Here ultraproduct
constructions provide a framework for comparing certain
test ideals in mixed characteristic with their positive
characteristic counterparts. In particular, we consider rings
of the form
$$R_i=W(k_i)[[X_2,\ldots,X_d]]
\quad\text{and}\quad
R^\sharp=W\bigl(\ulim k_i\bigr)[[X_2,\ldots,X_d]],$$
where the $k_i$ are perfect fields of a fixed positive
characteristic. For monomial ideals and fixed rational
exponents, we establish the comparison
$\bigl(\ulim\tau(R_i,\mathfrak a_i^t)\bigr)R^\sharp
=
\tau(R^\sharp,(\mathfrak a^\sharp)^t).$
The test ideals here are understood in the sense of
\cite[Definition 8.13]{BMP}. We also establish several
ultraproduct transfer results for trace-image descriptions,
Skoda-type identities, and colon containments.

Our next objective is to determine the Goto number more
precisely. The lower bound in Theorem~\ref{thm:WY} depends
only on the dimension and the order of the parameter ideal.
The associated graded ring, however, contains additional
information that can determine the exact value.

Recall that a standard graded Cohen--Macaulay algebra is called
level if its canonical module is generated in a single degree.
Suppose that
$\operatorname{gr}_{\mathfrak m}(R)$
is a standard graded Cohen--Macaulay level algebra with Hilbert
series
$\frac{h(t)}{(1-t)^d},r:=\deg h(t).$
For a parameter ideal
$J=(f_1,\ldots,f_d),f_i\in\mathfrak m^s\setminus\mathfrak m^{s+1},$
whose initial forms constitute a homogeneous system of
parameters, we obtain the following exact formulas:
$(J:\mathfrak m^t)=J+\mathfrak m^{r+d(s-1)+1-t},
\overline J=\overline{\mathfrak m^s},$
and, consequently,
${g(J)=(d-1)(s-1)+r.}$

Thus, the Goto number is governed not only by the dimension
and the order of the parameter ideal, but also by the degree
of the $h$-polynomial of the tangent cone. In particular, the
difference between the general singular lower bound and the
actual Goto number can be arbitrarily large.

We give examples illustrating this phenomenon, including
singular mixed-characteristic Gorenstein and non-Gorenstein
rings. In the regular case, the same computation recovers
$\overline Q=\mathfrak n^s,g(Q)=(d-1)(s-1),$
for
$Q:=(u_1^s,\ldots,u_d^s)$
in a regular local ring.
Then, we show that the Cohen--Macaulay hypothesis in
Question~\ref{q:nonCM} cannot be omitted. 
The final section  considers further questions from
\cite{WY}
concerning the interaction between test ideals, colon ideals,
and the associated graded ring, see  {Questions}    \ref{eq:q58} and \ref{eq:q86}. In particular, we investigate
the containments
$ J+\tau(R,J^{d-1})\subseteq (J:\overline J)$ and
$ \tau(R,\mathfrak m^t)
\subseteq \mathfrak m^{t+a(\operatorname{gr}_{\mathfrak m}(R))+1},$
in the settings specified below.

 \section{ Iterated socles under ultraproducts}
 
  Rings in this section are containing a field.
 
 \subsection{The Lefschetz approximation setup}

 \begin{definition}
 	Let $W$ be a set. An \emph{ultrafilter} on $W$ is a collection
 	$\U\subseteq\mathcal P(W)$ of subsets of $W$ such that:
 	\begin{enumerate}
 		\item $\varnothing\notin\U$ and $W\in\U$;
 		\item if $A\in\U$ and $A\subseteq B\subseteq W$, then $B\in\U$;
 		\item if $A,B\in\U$, then $A\cap B\in\U$;
 		\item for every $A\subseteq W$, exactly one of $A\in\U$ or
 		$W\setminus A\in\U$ holds.
 	\end{enumerate}
 	We say that a property holds \emph{for almost all $w$} (or
 	\emph{$\U$-almost everywhere}) if the set of $w$ for which it holds
 	belongs to $\U$.
 \end{definition}
 
 \begin{definition}
 	Let $(X_w)_{w\in W}$ be a family of nonempty sets. On the product
 	$\prod_{w\in W}X_w$ define
 	$(x_w)\sim(y_w)
 	\iff
 	\{w\in W : x_w=y_w\}\in\U .$
 	This is an equivalence relation. The \emph{ultraproduct} is
 	$	\ulim_{w} X_w
 	\;:=\;
 	\Bigl(\prod_{w\in W}X_w\Bigr)\Big/\sim .$
 	The class of $(x_w)$ is denoted $[x_w]$ or $\ulim_w x_w$. If all
 	$X_w=X$, one writes $X^{\U}$ and speaks of the \emph{ultrapower}.
 \end{definition}
 
 \begin{definition}
 	If each $R_w$ is a ring, the componentwise operations on
 	$\prod_w R_w$ descend to the ultraproduct
 	$ D(R)\;=\;\ulim_{w}R_w, $
 	which is again a ring with
 	$$[x_w]+[y_w]:=[x_w+y_w],\qquad
 	[x_w]\cdot[y_w]:=[x_w\cdot y_w],\qquad
 	0:=[0_w],\quad 1:=[1_w].	$$
 	More generally the same construction applies to groups, modules,
 	fields and first-order structures in a fixed language.
 \end{definition}
 
 \begin{fact}[{\L}o\'s's theorem]
 	A first-order sentence in the language of rings holds in
 	$\ulim_w R_w$ if and only if it holds in $R_w$ for almost all $w$.
 \end{fact}

 Let $(R,\m)$ be a local ring of equicharacteristic zero. Aschenbrenner--Schoutens
 \cite[Theorem~4.3, Theorem~5.2, Lemmas~5.3--5.4]{AS} construct a
 faithfully flat extension
 $R\to D(R)=\ulim_{w}R_w,$ 
 an ultraproduct of complete local rings $R_w$ of positive
 characteristic. Here $w$ runs over an index set $W$ equipped with an
 ultrafilter $\U$; Each $R_w$ is a complete local ring with maximal
 ideal $\m_w$, and the maximal ideal $\m$ of $R$ is approximated by the
 $\m_w$: the image of $\m$ in $D(R)$ is $\ulim_w \m_w$.
 The transfer principle says: given a \emph{fixed finite} list of
 elements $a_1,\dots,a_s\in R$ and ideals $I_1,\dots,I_t\subseteq R$,
 one can choose compatible approximations $a_{i,w}\in R_w$ and
 $I_{j,w}\subseteq R_w$ such that, for almost all $w$,
 \begin{itemize}
 	\item $\dim R_w=\dim R$ and $\operatorname{embdim}R_w=\operatorname{embdim}R$
 	(when finite);
 	\item $R_w$ is Cohen--Macaulay iff $R$ is;
 	\item a fixed tuple approximating a parameter sequence of $R$ is a
 	parameter sequence of $R_w$;
 	\item all fixed ideal relations among the $a_i$ and $I_j$  hold in $R_w$.
 \end{itemize}
 Let $x \in R$  and $I\subseteq R$ an ideal. We want to
 transfer statements such as:
 \begin{itemize}

 	\item $x$ is integral over an ideal $I$, i.e.\ $x$ satisfies an
 	equation $x^d+c_1x^{d-1}+\cdots+c_d=0$ with $c_i\in I^i$;
 	\item an element lies in the integral closure $\overline{I}$ of $I$.
 \end{itemize}

 \subsection{Dealing with {Question} \ref{q:charfree}}
 We give the additional numerical argument needed to transfer
 integral dependence.  In particular, we do not assume a uniform
 degree for integral equations in the approximating rings.
 
 \begin{lemma}\label{lem:length}
 	If $H$ is an $\m$-primary ideal of $R$, then, for almost all $w$,
 	$\length_{R_w}(R_w/H_w)=\length_R(R/H).$
 \end{lemma}
 \begin{proof}
 	Let the length on the right be $h$, and lift a composition series
 	of $R/H$ to ideals
 	$	H=H_0\subsetneq H_1\subsetneq\cdots\subsetneq H_h=R.$
 	Each factor is cyclic, nonzero, and annihilated by $\m$.
 	Choose its cyclic generator and finite generating sets of the
 	$H_i$.  The corresponding generation relations, inclusions
 	$\m H_i\subseteq H_{i-1}$, and nonmembership of the cyclic generator
 	in $H_{i-1}$ transfer to almost all $R_w$.
 	Nonmembership is preserved because the extension to $D(R)$ is
 	faithfully flat.  Since $\m_w$ is maximal, the resulting chain
 	has exactly $h$ simple factors.
 \end{proof}
 
 \begin{lemma}\label{lem:multiplicity}
 	Let $R$ be a complete Cohen--Macaulay local ring of
 	equicharacteristic zero, and let $I$ be an $\m$-primary ideal.  Then
 	$ \e(I_w,R_w)=\e(I,R) $
 	for almost all $w$.
 \end{lemma}
 \begin{proof}
 	The residue field of $R$ is infinite.  Choose a $d$-generated
 	minimal reduction $K\subseteq I$, where $d=\dim R$.  Thus $K$ is
 	a parameter ideal and, for some fixed integer $N\geq0$,
 	$ I^{N+1}=KI^N. $ 
 	This equality is a finite collection of ideal-membership
 	relations.  Consequently $K_w$ is a reduction of $I_w$ for
 	almost all $w$.  For these $w$, $R_w$ is Cohen--Macaulay of
 	dimension $d$ and $K_w$ is a parameter ideal.  Hence
 	$$	\e(I_w,R_w)
 	=\e(K_w,R_w)
 	=\length_{R_w}(R_w/K_w)
 	=\length_R(R/K)
 	=\e(K,R)
 	=\e(I,R).	$$
 	The middle equality is Lemma~\ref{lem:length}.
 \end{proof}
 
 \begin{lemma}
 	\label{lem:integraltransfer}
 	Let $R$ be as in Lemma~\ref{lem:multiplicity}.  Let
 	$J\subseteq I$ be $\m$-primary ideals.  If
 	$I_w\subseteq\ol{J_w}$ for almost all $w$, then
 	$I\subseteq\ol J$.
 \end{lemma}
 \begin{proof}
 	Integral dependence gives $\e(I_w,R_w)=\e(J_w,R_w)$.
 	Apply Lemma~\ref{lem:multiplicity} to $I$ and $J$ to obtain
 	$\e(I,R)=\e(J,R)$.  A complete Cohen--Macaulay local ring is
 	formally equidimensional.  Rees's multiplicity theorem
 	\cite[Theorem 11.3.1]{HSbook} now gives $I\subseteq\ol J$.
 \end{proof}
 
 \begin{theorem}\label{thm:field}
 	Let $(R,\m)$ be a Cohen--Macaulay local ring containing a field,
 	with $d=\dim R\geq2$.  If $s\geq2$ and $J\subseteq\m^s$ is a
 	parameter ideal, then
 	$ (J:\m^{(d-1)(s-1)})\subseteq\ol J. $
 	If $R$ is not regular, then also
 	$ (J:\m^{(d-1)(s-1)+1})\subseteq\ol J. $
 \end{theorem}
 \begin{proof}
 	In positive characteristic these are precisely \cite[Theorem \ref{thm:WY}]{WY}.
 	Suppose that $R$ has characteristic zero.
  	First pass to the $\m$-adic completion.  This preserves dimension,
 	Cohen--Macaulayness, and regularity.  It also preserves the
 	colon condition, and integral dependence descends by faithful
 	flatness.  For completeness, if $z$ becomes integral over
 	$J\widehat R$ and $I=J+(z)$, then for some $n$
 	$ 	I^{n+1}\widehat R=JI^n\widehat R. $
 	Faithful flatness gives $I^{n+1}=JI^n$, and hence $z\in\ol J$.
 	We may therefore assume that $R$ is complete.
 	Take $t=t_0$, or $t=t_0+1$ in the nonregular case, and fix
 	$z\in J:\m^t$.  Choose compatible approximations of $R$, $J$,
 	$\m$, and $z$.  For almost all $w$,
 	$$	\dim R_w=d,\qquad
 	J_w\text{ is a parameter ideal},\qquad
 	J_w\subseteq\m_w^s,\qquad
 	z_w\m_w^t\subseteq J_w,	$$
 	and $R_w$ is Cohen--Macaulay.
 	In the nonregular case, $\edim R>d$; dimension and embedding
 	dimension transfer, so that almost all $R_w$ are nonregular.
 	The positive-characteristic theorem yields
 	$z_w\in\ol{J_w}$ for almost all $w$.
 	In particular $z_w\in\m_w$ for almost all $w$, so $z\in\m$.
 	Thus $J+(z)$ is a proper $\m$-primary ideal.
 	Apply Lemma~\ref{lem:integraltransfer} to $J\subseteq J+(z)$.
 	It follows that $z\in\ol J$, as required.
 \end{proof}

 \section{Mixed char. }
 
  \subsection{Low dimensional approach}
 \begin{proposition} \label{thm:dimtwo}
 	Let $(R,\m)$ be a two-dimensional Cohen--Macaulay local ring.
 	For every parameter ideal $J\subseteq\m^s$, $s\geq2$,
 	$(J:\m^{s-1})\subseteq\ol J.$
 	If $R$ is not regular, then
 	$ (J:\m^s)\subseteq\ol J. 	$
 \end{proposition}
 
 \begin{proof}
 	If $R$ is regular, Heinzer--Swanson prove, without a
 	characteristic restriction,
 	$g(J)=\ord_\m(J)-1$, see
 	\cite[Theorem 2.2]{HS}.
 	Since $\ord_\m(J)\geq s$, this proves the first assertion.
 	If $R$ is nonregular, Wang's theorem
 	\cite[Theorems 1.1--1.2]{Wang}, also reproduced as
 	\cite[Theorem 1.1]{WY}, gives
 	$(J:\m^s)^2=J(J:\m^s).$
 	Thus $(J:\m^s)$ is integral over $J$.
 	The inclusion $(J:\m^{s-1})\subseteq (J:\m^s)$ gives the other
 	assertion.
 \end{proof}
 
 \begin{proposition} 
 	\label{prop:wangrange}
 	Let $R$ be Cohen--Macaulay and either nonregular of dimension
 	at least two, or regular of dimension at least three.
 	Let $J$ be a parameter ideal and put
 	$r=\ord_\m(J)\geq2$.  Then, in any characteristic,
 	$	g(J)\geq r.$
 \end{proposition}
 \begin{proof}
 	Apply Wang's theorem with exponent $r$ and use the increasing
 	chain $(J:\m^a)\subseteq (J:\m^b)$ for $a\leq b$.
 	For $d=3$, $s=2$, one has $t_0=2\leq r$.
 \end{proof}
 
 \subsection{The containment $J+\tau(R,J^{d-1})\subseteq( J:\cl J)$} 

\begin{notation} Fix an ultrafilter $\U$ on an index set and let $D:=\ulim R_i$.
 If $I_i\subseteq R_i$ are ideals, write
 $I_{\U}:=\ulim I_i
 :=\{[z_i]\in D:z_i\in I_i\text{ for $\U$-almost all }i\}.$
 This is an ideal of $D$, called the internal ideal of the family.
 It need not be finitely generated.\end{notation}
 
 \begin{lemma}\label{lem:internal}
 	Suppose $I_i=(f_{1,i},\ldots,f_{r,i})$ with the same finite $r$
 	for every $i$, and put $f_j=[f_{j,i}]$. Then
 $I_{\U}=(f_1,\ldots,f_r)D.$
 	For any other family $H_i$ of ideals, one has
 	$$	\ulim(I_iH_i)=I_{\U}H_{\U},\qquad
 	\ulim(I_i:H_i)=(I_{\U}:H_{\U}). 	$$
 	The equality of colons holds even without the bound on the
 	number of generators of $I_i$.
 \end{lemma}
 
 \begin{proof}
 	For the first assertion, express each element of $I_i$ as a sum
 	of the same $r$ generators with coefficients in $R_i$, and take
 	their ultraproducts. For the product assertion, every element of
 	$I_iH_i$ is $\sum_{j=1}^r f_{j,i}h_{j,i}$ with $h_{j,i}\in H_i$.
 	The fixed bound on the number of summands permits passage to $D$.
 	For the colon assertion, one inclusion is immediate. Conversely,
 	if $z_i\notin(I_i:H_i)$ on a set in $\U$, choose $h_i\in H_i$
 	with $z_ih_i\notin I_i$ on that set. Then $[h_i]\in H_{\U}$
 	and $[z_i][h_i]\notin I_{\U}$. This proves the other inclusion.
 \end{proof}
 
 \begin{definition}
 	For a family on which the mixed-characteristic test ideals are
 	defined, set
$T_{\U}(\ideal^t):=\ulim\tau(R_i,\ideal_i^t).$
 
 \end{definition}
 
 \begin{proposition}
 	\label{prop:trace}
 	Suppose that, for each $i$, a specified $R_i$-linear map
 	$ 	\phi_i:M_i\to R_i$
  satisfies 
$\operatorname{im}(\phi_i)=\tau(R_i,\ideal_i^t). $
 	Then the coordinatewise map
 $\phi_{\U}:\ulim M_i\to D$ defined  by 
 $[m_i]\mapsto[\phi_i(m_i)]$
 	is $D$-linear and satisfies
$T_{\U}(\ideal^t)=\operatorname{im}(\phi_{\U}). $
 \end{proposition}
 
 \begin{proof}
 	The map is well-defined and linear by its coordinatewise
 	definition. If $[z_i]\in T_{\U}(\ideal^t)$, choose
 	$m_i\in M_i$ with $\phi_i(m_i)=z_i$ on the set in $\U$
 	where $z_i\in\tau(R_i,\ideal_i^t)$. These choices give a preimage
 	of $[z_i]$. The reverse containment is immediate.
 \end{proof}

 \begin{proposition}\label{prop:skoda}
 	Let the $R_i$ be complete  normal $\mathbb Q$-Gorenstein
 	local domains of mixed characteristic. Suppose $\ideal_i$ is
 	generated by a fixed number $r$ of specified elements. Then, for
 	every rational $t\geq r$,
 	$T_{\U}(\ideal^t)=\ideal_{\U}\,T_{\U}(\ideal^{t-1}).$
 \end{proposition}
 
 \begin{proof}
 	Apply Skoda's theorem in each factor
 	\cite[Theorem 8.7 and Corollary 8.20(b)]{BMP} and then apply
 	Lemma~\ref{lem:internal} to the products. 
 \end{proof}
 
 \begin{proposition}\label{prop:one-sided}
 	Let $R$ be a complete  normal $\mathbb Q$-Gorenstein
 	local domain of mixed characteristic and dimension $d\geq2$.
 	For every parameter ideal $J$,  
 	$J+\tau(R,J^{d-1})\subseteq( J:\cl J).$
 	Consequently, for a family of such rings of fixed dimension $d$
 	and parameter ideals $J_i$, let
 	$ J_{\U}:=\ulim J_i$ and $ C_{\U}:=\ulim\cl{J_i}.$
 	Then
$J_{\U}+T_{\U}(J^{d-1})
 	\subseteq (J_{\U}:C_{\U}).$
 \end{proposition}
 
 \begin{proof}
 	We use monotonicity, invariance under integral closure of the
 	ideal, and Skoda's theorem. The integral-closure invariance for
 	the unified ideal follows from
 	\cite[Proposition 8.6(c), Theorem 8.11, and Definition 8.13]{BMP},
 	using a common sufficiently small perturbation for the two
 	ideals in this single ring.
 	For any ideal $L$ and $t\geq0$,
$	L\tau(R,L^t)\subseteq\tau(R,L^{t+1}). $
 	Indeed, for $f\in L$, principal Skoda gives
 	$f\tau(R,L^t)=\tau(R,(f)^1L^t)$, and monotonicity gives the
 	containment in $\tau(R,L^{t+1})$. Thus
 	\begin{align*}
 		\cl J\,\tau(R,J^{d-1})
 	=\cl J\,\tau(R,(\cl J)^{d-1})
 	\subseteq\tau(R,(\cl J)^d)
 		=\tau(R,J^d)
 		=J\tau(R,J^{d-1})\subseteq J.
 	\end{align*}
 	Here Skoda applies to the $d$ generators of $J$. This proves the
 	single-ring assertion. Take ultraproducts and use the colon
 	equality in Lemma~\ref{lem:internal} for the second assertion.
 \end{proof}

 \subsection{An unconditional mixed-characteristic comparison}
 \label{sec:comparison}
 
 Fix a prime $p$, an integer $d\geq1$, and perfect fields $k_i$
 of characteristic $p$.  
 \begin{notation}
 Let $R_i:=W(k_i)[[X_2,\ldots,X_d]],
 \mm_i=(p,X_2,\ldots,X_d),$ and 
 $K:=\ulim k_i.$
Let $D:=\ulim R_i$ and $\mm:=(p,X_2,\ldots,X_d)D$.
 \end{notation}
 The symbols $X_j$ on $D$ denote the classes of the corresponding
 variables. The field $K$ is perfect. In the formulas below, set
 $u_1=p$ and $u_j=X_j$ for $2\leq j\leq d$.
 
 \begin{lemma}\label{lem:witt}
 	There is a natural isomorphism
 $R^\sharp:=\varprojlim_N D/\mm^N
 	\simeq W(K)[[X_2,\ldots,X_d]].$
 	It identifies the images of $p,X_2,\ldots,X_d$ with the indicated
 	regular parameters. In particular, $R^\sharp$ is a complete
 	Noetherian regular local ring of mixed characteristic $(0,p)$.
 \end{lemma}
 
 \begin{proof}
 	For every fixed $N$, Lemma~\ref{lem:internal} gives
 $D/\mm^N\simeq\ulim(R_i/\mm_i^N).$
 	Only finitely many monomials $X^\alpha$ with $|\alpha|<N$
 	are needed to describe $R_i/\mm_i^N$; their coefficients are
 	respectively taken modulo $p^{N-|\alpha|}$. Thus they are
 	truncated Witt vectors of bounded lengths.
For each fixed $n$, addition and multiplication on $W_n(k_i)$
 	are given by the same finite list of Witt polynomials, because
 	$p$ is fixed. Consequently
 $\ulim W_n(k_i)=W_n(K).$
 	Applying this simultaneously to the finitely many coefficient
 	slots gives compatible ring isomorphisms
 	$D/\mm^N\simeq
 	W(K)[[X_2,\ldots,X_d]]/(p,X_2,\ldots,X_d)^N.$
 	Take the inverse limit over $N$. This argument constructs the
 	map $D\to R^\sharp$; no approximation theorem asserting faithful
 	flatness of that map is needed or claimed.
 \end{proof}
 
 Let $E\subset\NN^d$ be a fixed finite nonempty set. Assume the
 monomial ideals
 $$\ideal_i=(u^\alpha:\alpha\in E)R_i,
 \quad \ideal_D=(u^\alpha:\alpha\in E)D,
 \quad \ideal^\sharp=(u^\alpha:\alpha\in E)R^\sharp$$
 are proper. Their common Newton polyhedron is
$ P=\conv(E)+\mathbb R_{\geq0}^{d}.$
 We record a version of Dickson's lemma that we need:
 \begin{lemma}[Dickson's lemma]
 	Let $d\geq 1$ and let $S\subseteq\mathbb{N}^d$.
 	Then $S$ has only finitely many minimal elements with respect
 	to the coordinatewise partial order. Moreover, if $S$ is upward
 	closed, then every element of $S$ lies above a minimal element
 	of $S$.
 \end{lemma}

 \begin{theorem}\label{thm:monomial}
 	For every   rational $t>0$,  
 $T_{\U}(\ideal^t) =\bigl(u^\beta:\beta\in\NN^d,
 \beta+\one\in\Int(tP)\bigr)D$ and 
$T_{\U}(\ideal^t)R^\sharp=\tau(R^\sharp,(\ideal^\sharp)^t).$
 	In particular, the first ideal is finitely generated by a set of  monomials independent of $i$.
 \end{theorem}
 
 \begin{proof}
 	The mixed-characteristic monomial formula
 	\cite[Proposition 8.21]{BMP} says that
 	$$
 	\tau(R_i,\ideal_i^t)
 	=\bigl(u^\beta:\beta+\one\in\Int(tP),\ \beta\in\NN^d\bigr)R_i.
 	$$
 	The exponent set is upward closed for the coordinatewise order.
 	Dickson's lemma gives a finite set $M\subset\NN^d$ of minimal
 	exponents generating it. Both $P$ and $M$ are independent of $i$.
 	The first assertion now follows from Lemma~\ref{lem:internal}.
 	Lemma~\ref{lem:witt} identifies $R^\sharp$ as another complete
 	regular local ring with the same regular parameters. Applying
 	the same monomial formula to this ring proves the second
 	assertion.
 \end{proof}

 \begin{corollary}
 	For fixed $E$ and $t>0$ as above, there exists $\epsilon_0>0$
 	such that, for every $i$ and rational $0<\epsilon<\epsilon_0$,
 $\tau(R_i,\ideal_i^t)=\tau(R_i,\ideal_i^{t+\epsilon}).$
 \end{corollary}
 
 \begin{proof}
 	Write the nontrivial supporting inequalities of $P$ as
 	$\langle v_j,x\rangle\geq c_j$, with integral vectors $v_j$ and
 	positive integers $c_j$; there are finitely many of them.
 	For $\beta\in\NN^d$, the numbers
 	$\langle v_j,\beta+\one\rangle$ are integers. Choose
 $	0<\epsilon_0<
 	\min_j\frac{\lfloor tc_j\rfloor+1-tc_j}{c_j}.$
 	Then replacing $t$ by $t+\epsilon$ changes none of the strict
 	integer inequalities defining the monomial test ideal.
 	Coordinate inequalities with zero right-hand side do not change.
 \end{proof}
 
 \subsection{The Watanabe--Yoshida identity for equal parameter powers}
 
We continue with the fixed-prime family in
 Section~\ref{sec:comparison}, and assume $d,s\geq2$. Set
 $$
 Q_i:=(p^s,X_2^s,\ldots,X_d^s)R_i,
 \quad Q:=(u_1^s,\ldots,u_d^s)D,
 \quad Q^\sharp=QR^\sharp,
 \quad h:=(d-1)(s-1).
 $$
 
 \begin{theorem}\label{thm:WY2}
 	With ordinary ideal integral closure on $D$, one has
 $	\cl Q=\mm^s=\ulim\cl{Q_i}$, 
$T_{\U}(Q^{d-1})=\mm^h$ and
${( Q:\cl Q)=Q+T_{\U}(Q^{d-1}).}$
Also, $
  {(Q^\sharp:\cl{Q^\sharp})
 		=Q^\sharp+\tau(R^\sharp,(Q^\sharp)^{d-1}).} $
 	Furthermore,
 	$ (	Q:\mm^h)=\mm^s=	\cl Q$, and $(Q:\mm^{h+1})=Q+\mm^{s-1}\not\subseteq\cl Q. $
 \end{theorem}
 
 \begin{proof}
 	We first record computations valid in every regular local ring
 	$(S,\mathfrak n)$ of dimension $d$ with parameters
 	$v_1,\ldots,v_d$, in any characteristic. Let
 	$Q_S:=(v_1^s,\ldots,v_d^s)$.
 	Every monomial $v^\alpha$ of total degree $s$ is integral over
 	$Q_S$, since
 $(v^\alpha)^s\in Q_S^s.$
 	Also $\mathfrak n^s$ is integrally closed. Indeed,
 	$\gr_{\mathfrak n}S$ is a polynomial domain, so the
 	$\mathfrak n$-adic order is multiplicative. An integral equation
 	over $\mathfrak n^s$ cannot hold for an element of order less
 	than $s$, because its highest power would be the unique term of
 	smallest order. Therefore
 $\cl{Q_S}=\mathfrak n^s.$
 The associated graded ring of $S/Q_S$ is
$k[V_1,\ldots,V_d]/(V_1^s,\ldots,V_d^s),$
 	an artinian graded Gorenstein ring with socle degree
 	$L=d(s-1)$. Its perfect pairings into degree $L$ give
$(Q_S:\mathfrak n^r)=Q_S+\mathfrak n^{L+1-r}$
 for $0\leq r\leq L+1$.
 	To justify the filtered step, an element of order $j\leq L-r$
 	has nonzero initial form which can be multiplied by a degree-$r$
 	form to give a nonzero product. It therefore does not annihilate
 	$\mathfrak n^r$ modulo $Q_S$. The opposite inclusion follows
 	from $\mathfrak n^{L+1}\subseteq Q_S$.
 Let us apply these computations in each $R_i$. The fixed finite
 	generating sets for powers of $\mm_i$ allow all the displayed
 	colon identities to pass to $D$ by Lemma~\ref{lem:internal}.
 	It remains to check that the internal integral closure is the
 	ordinary integral closure in $D$. If $z\in\cl Q$, choose a
 	single integral equation for $z$ over $Q$. Its degree is finite,
 	and membership of its coefficients in the fixed powers of $Q$
 	is witnessed by finitely many equations. It follows that
 	$z_i\in\cl{Q_i}=\mm_i^s$ almost everywhere, so $z\in\mm^s$.
 	Conversely, the degree-$s$ monomial equations above work directly
 	in $D$, and the integral closure of an ideal is an ideal. Hence
 	$\mm^s\subseteq\cl Q$.
 The Newton polyhedron of $Q_i$ is
 	$\{x\in\mathbb R_{\geq0}^d:x_1+\cdots+x_d\geq s\}.$
 	Theorem~\ref{thm:monomial} with $t=d-1$ shows that its test
 	ideal is generated by monomials with
 	$|\beta|+d>s(d-1)$, i.e., 
 	$|\beta|\geq h.$
 	Thus $T_{\U}(Q^{d-1})=\mm^h$. Setting $r=s$ in the colon
 	formula gives $( Q:\mm^s)=Q+\mm^h$, proving the first boxed
 	identity. The same calculations apply to $R^\sharp$, proving
 	the second one. Finally set $r=h$ and $r=h+1$ in that formula.
 	The element $u_1^{s-1}$ belongs to the latter colon and not to
 	$\mm^s$, in each factor and hence also in $D$.
 \end{proof}
 
 \begin{remark}
 We will present in Proposition 
 \ref{prop:regular}  another argument of {Theorem} \ref{thm:WY2}  without any use of mixed-characteristic
 	test ideal from \cite{BMP}.
 
 \end{remark}

 \begin{proposition}\label{prop:no-lefschetz}
 	Let $(R,\mm)$ have mixed characteristic $(0,p)$. There is no
 	faithfully flat unital map from $R$ to an ultraproduct of nonzero
 	rings of prime positive characteristics.
 \end{proposition}
 
 \begin{proof}
 	If the factors have characteristic $p$ on a set in $\U$, then
 	$p$ maps to zero, contradicting the injectivity of a faithfully
 	flat map. Otherwise $p$ is a unit in almost every factor, and
 	hence in the ultraproduct. This contradicts faithful flatness,
 	since the proper ideal $pR$ would extend to the unit ideal.
 \end{proof}

 \section{The role of Cohen--Macaulayness}
 
 We show   the basic containment of Theorem \ref{thm:WY} fails for $T$.
 The ring $T$ is reduced, complete, two-dimensional, and
 not Cohen--Macaulay; moreover, $\pi$ is a nonzerodivisor.
 Choosing $V=k[[\pi]]$ gives examples over every field.
 
 \begin{proposition}
 	\label{prop:noncm}
 	Let  $(V,\pi)$ be a complete discrete valuation ring.  Define
 	$ T:=V[[x,y]]/(xy,y^2-\pi y)$ and $ t=(\pi,x,y)$.
 	For every $s\geq3$, the ideal
 	$ J=(\pi^s,x^s)T $
 	is a parameter ideal contained in $\mathfrak t^s$, but
 	$ y^{s-1}\in (J:\mathfrak t^{\,s-1})\setminus \ol J.$
 	
 \end{proposition}
 \begin{proof}
 	In $V[[x,y]]$, one has
 	$ (xy,y^2-\pi y)=(y)\cap(x,y-\pi). $
 	Both ideals on the right are prime.  Their quotient rings
 	are $V[[x]]$ and $V$, of dimensions two and one.
 	Thus $T$ is reduced and has minimal components of different
 	dimensions.  Since it is complete, it cannot be
 	Cohen--Macaulay.  Neither minimal prime contains $\pi$,
 	so $\pi$ is a nonzerodivisor.  In mixed characteristic
 	$p$ is a unit times a positive power of $\pi$, so $T$ is
 	also $p$-torsionfree.
 	
 	The relations in $T$ give $y^n=\pi^{n-1}y$ for $n\geq1$.
 	In particular $y^{s+1}=\pi^sy\in J$.
 	Together with $\pi^s,x^s\in J$, this proves that
 	$J$ is $\mathfrak t$-primary and hence a parameter ideal.
 	Let $a:=y^{s-1}=\pi^{s-2}y$.  Any degree-$(s-1)$
 	monomial in $\pi,x,y$ containing $x$ is killed by $a$.
 	For a monomial $\pi^iy^j$ with $i+j=s-1$, one has
 	$ 	a\pi^iy^j=\pi^{2s-3}y\in(\pi^s)T\subseteq J, $
 	since $2s-3\geq s$ for $s\geq3$.
 	Thus $a\mathfrak t^{s-1}\subseteq J$.
 	
 	Under the quotient map
 	$ 	T\to T/(x,y-\pi)T\cong V $
 	the ideal $J$ maps to $(\pi^s)$ and $a$ maps to $\pi^{s-1}$.
 	An integral equation descends under any ring homomorphism.
 	But $(\pi^s)$ is integrally closed in the discrete valuation
 	ring $V$, and does not contain $\pi^{s-1}$.
 	Consequently $a\notin\ol J$.
 \end{proof}
 
 \section{Bounds on Goto number}
 
  We start with some lemmas.
  
 \begin{lemma} \label{lem:filtered}
 	Let $(A,\mathfrak q)$ be  artinian  so that
 	$B:=\gr_{\mathfrak q}(A)$ is a standard graded artinian algebra
 	whose socle is concentrated in degree $L$.
 	Then
 	$ (	0:_A\mathfrak q^t)=\mathfrak q^{L+1-t}$
 	for $0\leq t\leq L+1$.
 \end{lemma}
 \begin{proof}
 	Because $H_i=0$ for $i>L$, Nakayama's lemma gives
 	$\mathfrak q^{L+1}=0$.  Thus the right-hand side is contained
 	in the left-hand side.  For the reverse inclusion let
 	$0\neq a\in A$ with $j=\ord_{\mathfrak q}(a)\leq L-t$.
 	Its initial form $a^*\in B_j$ is nonzero.  Since the socle
 	of $B$ is concentrated in degree $L$, one can successively
 	multiply $a^*$ by degree-one elements without obtaining zero
 	until reaching degree $L$.  In particular, there is
 	$b^*\in B_t$ with $a^*b^*\neq0$.
 	Lift $b^*$ to $b\in\mathfrak q^t$.  The initial form of $ab$
 	is $a^*b^*\neq0$, so $ab\neq0$.
 	Hence $a\notin0:_A\mathfrak q^t$, as required.
 	The endpoint cases $t=0,L+1$ follow by the same argument.
 \end{proof}
 
 \begin{lemma}
 	\label{lem:powernotintegral}
 	Let $(R,\m)$ be a   local ring of positive dimension,
 	and let $s\geq2$.  Then
 	$ \m^{s-1}\not\subseteq\ol{\m^s}. $
 \end{lemma}
 \begin{proof}
 	Otherwise, the ideal $\m^s$ would be a reduction of
 	$\m^{s-1}$.  For sufficiently large $N$,
 	$ (\m^{s-1})^{N+1}=\m^s(\m^{s-1})^N, $
 	which says $\m^{(s-1)(N+1)}=\m^{(s-1)(N+1)+1}$.
 	Nakayama's lemma would force a positive power of $\m$
 	to vanish, contradicting $\dim R>0$.
 \end{proof}
 The following has application to {Question} \ref{q:charfree} when the ring is of mixed char.
 \begin{theorem}
 	\label{thm:level}
 	Let $(R,\m,k)$ be a $d$-dimensional Cohen--Macaulay local
 	ring, $d\geq2$.  Suppose  
 	$G:=\gr_\m(R)$ is a Cohen--Macaulay, standard graded
 	\emph{level} $k$-algebra, with Hilbert series
 	$	\operatorname{Hilb}_G(t)=\frac{h(t)}{(1-t)^d}$ with $ r:=\deg h(t).$
 	Fix $s\geq2$ and a parameter ideal
 	$J=(f_1,\ldots,f_d)$, where each
 	$f_i\in\m^s\setminus\m^{s+1}$ and the initial forms
 	$f_1^*,\ldots,f_d^*\in G_s$ constitute a homogeneous
 	system of parameters in $G$.  Then for
 	$0\leq t\leq r+d(s-1)+1$,
 	
 	\begin{enumerate}[label=\textup{(\roman*)}]
 		\item  
 		$(	J:\m^t)=J+\m^{r+d(s-1)+1-t}$,
 		\item
 		$	\ol J=\ol{\m^s}$ and
 		$g(J)=(d-1)(s-1)+r.$ 
 	\end{enumerate}	
 \end{theorem}
 \begin{proof}
 	The initial forms form a regular sequence in the
 	Cohen--Macaulay ring $G$.  The standard lifting theorem for
 	regular initial sequences gives
 	$ \gr_\m(R/J)\cong G/(f_1^*,\ldots,f_d^*). $
 	This Artinian graded ring is level: for a Cohen--Macaulay
 	level graded ring, quotienting by a homogeneous regular
 	sequence preserves levelness.  One way to see this is to
 	use the graded canonical module: it is generated in a single
 	degree, and the canonical module of the quotient is the
 	quotient of the canonical module, up to a grading shift.
 	The Hilbert series of the quotient is
 	$ 	h(t)(1+t+\cdots+t^{s-1})^d, $
 	so its top socle degree is $L=r+d(s-1)$.
 	Lemma~\ref{lem:filtered}, applied to $R/J$, proves
 	(i).
 	
 	We next verify that $J$ is a reduction of $\m^s$.
 	Since $G/(f_1^*,\ldots,f_d^*)$ has finite length,
 	for all sufficiently large $n$ one has
 	$G_n=\sum_i f_i^*G_{n-s}$.  Lifting to $R$ gives
 	$ \m^n=J\m^{n-s}+\m^{n+1}. $
 	Nakayama's lemma implies $\m^n=J\m^{n-s}$.
 	Taking $n=qs$ for large $q$, we obtain
 	$(\m^s)^q=J(\m^s)^{q-1}$, as required.
 	Consequently $\ol J=\ol{\m^s}$.
 	Let  $t_*:=r+(d-1)(s-1)$.  Formula~(i)
 	gives $(J:\m^{t_*})=J+\m^s\subseteq\ol J$.
 	At the next exponent the same formula gives
 	$(J:\m^{t_*+1})=J+\m^{s-1}$, which is not contained
 	in $\ol J$ by Lemma~\ref{lem:powernotintegral}.
 	Thus $g(J)=t_*$.
 	Finally $r\geq0$.  If $R$ is nonregular, then $r\geq1$:
 	indeed, $r=0$ would imply $e(R)=h(1)=1$, forcing the
 	Cohen--Macaulay local ring $R$ to be regular.
 \end{proof}
 
 \begin{corollary}
 	\label{cor:generic}
 	Under the tangent-cone hypotheses of Theorem~\ref{thm:level},
 	if the residue field $k$ is infinite, the $d$-tuples
 	$(f_1^*,\ldots,f_d^*)\in(G_s)^d$ that satisfy the theorem
 	form a nonempty Zariski-open subset.  Thus the exact formula
 	\ref{thm:level}(ii) holds for a general $d$-generated
 	parameter ideal of order $s$.
 \end{corollary}
 \begin{proof}
 	A standard graded Cohen--Macaulay algebra over an infinite
 	field admits a linear homogeneous system of parameters
 	$\ell_1,\ldots,\ell_d$.  Their $s$-th powers are a
 	homogeneous system of parameters of degree $s$, so the
 	asserted locus is nonempty.  The homogeneous-system-of-
 	parameters condition is open in $(G_s)^d$ (equivalently,
 	the condition that the quotient have dimension zero is open).
 	Theorem~\ref{thm:level} applies to every lift of such a tuple.
 \end{proof}

 \begin{proposition} 
 	\label{prop:regular}
 	Suppose $(V,\pi)$ is a  complete mixed-characteristic discrete
 	valuation ring. 	Let $(A,\n)$ be any regular local ring of dimension $d\geq2$,
 	with regular system of parameters $u_1:=\pi,\ldots,u_d$.
 	For $s\geq2$, let
 	$	Q:=(u_1^s,\ldots,u_d^s).$ 
 	Then
 	$\ol Q=\n^s$ and
 	$g(Q)=(d-1)(s-1). $
 	
 \end{proposition}
 \begin{proof}
 	The associated graded ring of $A/Q$ is
 	$ k[U_1,\ldots,U_d]/(U_1^s,\ldots,U_d^s), $
 	an Artinian graded complete intersection of socle degree
 	$L=d(s-1)$.  The local ring $A/Q$ is Gorenstein.
 	Lemma~\ref{lem:filtered} gives
 	$
 	(Q:\n^t)=Q+\n^{d(s-1)+1-t}$ for all $0\leq t\leq d(s-1)+1.$
 	Every degree-$s$ parameter monomial $u^\alpha$ is integral over
 	$Q$, since $(u^\alpha)^s\in Q^s$.  Hence $\n^s$ is integral
 	over $Q$.
 	Also $\n^s$ is integrally closed: the $\n$-adic order is
 	multiplicative because $\gr_\n A$ is a polynomial domain,
 	and taking orders in an integral equation excludes elements
 	of order less than $s$.  Thus $\ol Q=\n^s$.
 	For $t_0=(d-1)(s-1)$ the displayed colon is $\n^s$.
 	For $t=t_0+1$ it is $Q+\n^{s-1}$ and contains
 	$u_1^{s-1}\notin\ol Q$.
 	This proves the exact value of $g(Q)$.
 \end{proof}
 Theorem \ref{thm:WY} holds for these parameter ideals in every
 residue characteristic.  For $e=2$, the exponent
 $(d-1)(s-1)+1$ is best possible even among singular
 mixed-characteristic Cohen--Macaulay domains:
 \begin{example}
 	\label{hypersurface}
 	Let $A$ be either $k[[u_1,\ldots,u_d]]$ over any field,
 	or $V[[u_2,\ldots,u_d]]$ over a complete mixed-characteristic
 	discrete valuation ring $(V,\pi)$, with $u_1=\pi$.  Fix $d,s,e\geq2$.
 	Set
 	$(R,\m):=A[[z]]/(z^e-u_1^{e-1}u_2)$ and 
 	$J:=(u_1^s,\ldots,u_d^s)R. $
 	Then $R$ is a singular complete Cohen--Macaulay  domain,
 	$\dim R=d$,  
 	$ \ol J=\m^s$ and
 	$	g(J)=(d-1)(s-1)+e-1. $
 	
 \end{example}
 \begin{proof}
 	The ring is a hypersurface of dimension $d$ and embedding
 	dimension $d+1$, so it is Cohen--Macaulay and nonregular.
 	Its associated graded ring is
 	$$
 	G:=\gr_\m (R)
 	=k[U_1,\ldots,U_d,Z]/(Z^e-U_1^{e-1}U_2).
 	$$
 	The displayed binomial is irreducible in every characteristic:
 	view it as a primitive polynomial of degree one in $U_2$
 	over the unique factorization domain
 	$k[U_1,U_3,\ldots,U_d,Z]$.
 	Consequently, $G$ is a domain.  Separatedness of the
 	$\m$-adic filtration shows that $R$ is a domain, and its
 	$\m$-adic order is multiplicative.
 	The ring $G$ is free over $k[U_1,\ldots,U_d]$, with basis
 	$\{1,Z,\ldots,Z^{e-1}\}$.  Thus the initial forms
 	$U_1^s,\ldots,U_d^s$ form a regular sequence in $G$, and
 $$\gr_\m(R/J)
 	=k[U_1,\ldots,U_d,Z]/
 	(U_1^s,\ldots,U_d^s,Z^e-U_1^{e-1}U_2).	$$
 	This is an artinian graded complete intersection with Hilbert
 	series
 	$ (1+t+\cdots+t^{s-1})^d
 	(1+t+\cdots+t^{e-1}) $
 	and socle degree $L=d(s-1)+e-1$.
 	The local ring $R/J$ is also artinian Gorenstein.
 	It follows from Lemma~\ref{lem:filtered} that
 	$$(\ast):\quad
 ( J:\m^t)=J+\m^{d(s-1)+e-t}
 	\qquad(0\leq t\leq d(s-1)+e).
 	$$
 	
 	Write $\n=(u_1,\ldots,u_d)A$.  The equation for $z$ shows
 	that $z$ is integral over $\n R$, so $\m$ is integral over
 	$\n R$.  Also $\n^sR$ is integral over $J$, by the same
 	monomial equations as in Proposition~\ref{prop:regular}.
 	Powers and transitivity of integral dependence therefore
 	show that $\m^s$ is integral over $J$.
 	Since $G$ is a domain, the order argument used above shows
 	that $\m^s$ is integrally closed.  Hence $\ol J=\m^s$.
 	For $g=(d-1)(s-1)+e-1$, formula $(\ast)$ gives
 	$( J:\m^g)=\m^s$.  At $g+1$ it gives $J+\m^{s-1}$,
 	which contains $u_1^{s-1}\notin\ol J$.
 	The claimed Goto number follows.
 \end{proof}

 The singular bound in Theorem~\ref{thm:WY} is sharp in
 mixed characteristic even for non-Gorenstein domains,
 in every dimension $d\geq2$:
 
 \begin{example}
 	\label{determinantal}
 	Let $(V,\pi)$ be a complete discrete valuation ring of mixed
 	characteristic, let
 	$A:=V[[x_2,\ldots,x_d]]$ and
 	$(R, \m):=\frac{A[[u,v]]}
 	{(u^2-\pi v,\ uv-\pi x_2,\ v^2-u x_2)}$ and	$J:=(\pi^s,x_2^s,\ldots,x_d^s)R.$
 	Then $R$ is a complete, Cohen--Macaulay  domain of
 	dimension $d$, but $R$ is  not Gorenstein.
 	For these parameter ideals,
 	$ \ol J=\m^s$ and $	g(J)=(d-1)(s-1)+1.$
  \end{example}
 \begin{proof}
 	The three relations reduce all monomials in $u,v$ to
 	$A$-linear combinations of $1,u,v$.  To prove independence,
 	consider the finite $A$-algebra
 	$	S=A[a,b]/(a^3-\pi,\ b^3-x_2).$
 	It is a complete regular local ring of dimension $d$,
 	with regular parameters $a,b,x_3,\ldots,x_d$; in particular
 	it is a domain.  As an $A$-module it is free with basis
 	$a^i b^j$ for $0\leq i,j<3$.  The substitutions
 	$u\mapsto a^2b$, $v\mapsto ab^2$ respect the displayed
 	relations, and $1,a^2b,ab^2$ are $A$-linearly independent.
 	Consequently $R$ is a domain and is free over the regular
 	local ring $A$ with basis $1,u,v$.
 	In particular $R$ is Cohen--Macaulay of dimension $d$.
 	The relations are homogeneous of degree two for the
 	maximal-ideal filtration, and the displayed free basis
 	is compatible with that filtration.  Writing
 	$P,X_2,\ldots,X_d,U,V$ for initial forms, we obtain
 	$$G=\gr_\m ( R)\cong
 	\frac{k[P,X_2,\ldots,X_d,U,V]}
 	{(U^2-PV,\ UV-PX_2,\ V^2-UX_2)}.$$
 	It is the homogeneous coordinate ring of the rational
 	normal cubic curve, with the extra polynomial variables
 	$X_3,\ldots,X_d$.  In particular $G$ is a standard
 	Gorenstein-free (non-Gorenstein), Cohen--Macaulay level
 	domain, and as a module over
 	$k[P,X_2,\ldots,X_d]$ it has homogeneous basis $1,U,V$
 	of degrees $0,1,1$.  More explicitly, the quotient of
 	$G$ by the linear regular sequence
 	$P,X_2,\ldots,X_d$ is
 	$k[U,V]/(U,V)^2$, with socle concentrated in degree one;
 	this verifies levelness.  Thus
 	$ \operatorname{Hilb}_G(t)=\frac{1+2t}{(1-t)^d},$ and so $r=1$.
 	The initial forms $P^s,X_2^s,\ldots,X_d^s$ are a homogeneous
 	system of parameters, so Theorem~\ref{thm:level} gives
 	$g(J)=(d-1)(s-1)+1$ and $\ol J=\ol{\m^s}$.
 	As $G$ is a domain, the $\m$-adic order on $R$ is
 	multiplicative, so each $\m^s$ is integrally closed;
 	hence $\ol J=\m^s$.
 	Finally, $R/(\pi,x_2,\ldots,x_d)R$ is
 	$k[u,v]/(u,v)^2$, whose socle has dimension two.
 	Therefore $R$ has Cohen--Macaulay type two and is not
 	Gorenstein.  Its embedding dimension is $d+2>d$, so it is
 	also nonregular.
 \end{proof}
 \section{Additional investigations}
 
 \subsection{A non containment $\tau(R,\mm^t)\nsubseteq\mm^{t+a(\gr_{\mm}R)+1}$ }
 \label{sec:q58}
 
 In this subsection, $\tau$ denotes the classical positive-characteristic
 Hara--Yoshida test ideal used in \cite{WY}. They asked:

 \begin{question}\label{eq:q58} (See \cite[Question 5.8]{WY}). 	When
 	$\tau(R,\mm^t)\subseteq\mm^{t+a(\gr_{\mm}R)+1}?$ 
 \end{question}
 
 Normality is not assumed
 in the  question.  
 The following example disproves the
 containment, even for a complete Gorenstein domain
 with Cohen--Macaulay associated graded ring.

 \begin{proposition}\label{thm:q58-counter}
 	Let $k=\overline{\mathbb F_p}$ for any prime $p$, and set
 	$(R,\mm):=k[[u^5,u^6,u^9]]$.
Then $R$ is a one-dimensional complete Gorenstein domain,
 	$G:=\gr_{\mm}(R)$ is Cohen--Macaulay, and
  $\operatorname{Hilb}_G(z)=\frac{1+2z+z^2+z^3}{1-z},
 a(G)=2.$
 	For every integer $t\geq0$,
 	$\tau(R,\mm^t)=u^{5t+14}k[[u]],$ and 
 	$u^{5t+14}\notin\mm^{t+3}.	$ 
 	In particular,
${u^{19}\in\tau(R,\mm)\setminus\mm^4.}$
\end{proposition}
 
 \begin{proof}Let $S:=k[[u]]$.
Let $x:=u^5$, $y:=u^6$, and $z:=u^9$. The numerical semigroup $H:=\langle5,6,9\rangle$. Recall that if $H$ is a numerical semigroup and $n\in H\setminus\{0\}$, the Ap\'ery set of $H$ with respect to $n$ is defined by $\operatorname{Ap}(H,n):=\{h\in H\mid h-n\notin H\}.$
Equivalently, it consists of the least elements of $H$ in each residue class modulo $n$. Since $H=\langle5,6,9\rangle$, we have
$\operatorname{Ap}(H,5)=\{0,6,12,18,9\},$
where the displayed elements are the least elements of $H$ in the residue classes $0,1,2,3,4$, respectively.
In residue classes $0,1,2,3,4$ modulo $5$, these are the least
	members of $H$. Consequently $R$ is free of rank $5$ over
$k[[x]]$, with basis $1,y,y^2,y^3,z$.
There is a presentation
$	R\simeq k[[X,Y,Z]]/(X^3-YZ,\ Z^2-Y^3).$
Here is a verification that this presentation has no extra
component. The quotient on the right has dimension one and
is a complete intersection: modulo $X$ its quotient is
$k[[Y,Z]]/(YZ,Z^2-Y^3)$, of length $5$.
It is therefore Cohen--Macaulay, $X$ is a nonzerodivisor, and
it is free of rank $5$ over $k[[X]]$. Its natural surjection
onto $R$, also free of rank $5$, is an isomorphism. In particular,
$R$ is Gorenstein.
For each $v\in H$, let
$\operatorname{ord}_{\mathfrak m}(u^v)
=\max\{r\geq 0\mid u^v\in\mathfrak m^r\}.$
Recall that
$\operatorname{Ap}(H,5)=\{0,6,12,18,9\}.$
Every $v\in H$ has a unique expression
$$
v=5q+w,\qquad q\geq 0,\quad
w\in\operatorname{Ap}(H,5).
$$
The $\mathfrak m$-adic orders of the corresponding Apéry monomials are
$$
\begin{array}{c|ccccc}
	w&0&6&12&18&9\\ \hline
	\operatorname{ord}_{\mathfrak m}(u^w)&0&1&2&2&1.
\end{array}
$$
Indeed, these orders follow directly from
$$u^6=y,\quad u^{12}=y^2,\quad
u^{18}=z^2,\quad u^9=z,$$
and the fact that $u^{18}\notin\mathfrak m^3$.
We claim that
$\operatorname{ord}_{\mathfrak m}(u^{5q+w})
=q+\operatorname{ord}_{\mathfrak m}(u^w).$
To see this, multiplication by $x=u^5$ increases the order by at least one. Conversely, the valuation of every element of $\mathfrak m^r$ is at least $5r$, and the displayed Ap\'ery orders give the required upper bound. Thus equality follows.
Every element $f\in R$ has a unique expression
$$f=\sum_{w\in\operatorname{Ap}(H,5)}u^w f_w(x),
\qquad f_w(x)\in k[[x]].$$
Since the five summands have distinct valuations modulo $5$, there is no cancellation between their terms of lowest valuation. It follows that
$\operatorname{ord}_{\mathfrak m}(xf)
=\operatorname{ord}_{\mathfrak m}(f)+1$
for every nonzero $f\in R$.
Consequently, multiplication by $x^*$ is injective on
$G=\operatorname{gr}_{\mathfrak m}(R).$
Since $\dim G=1$, we conclude that $G$ is Cohen--Macaulay.
Since $x^*$ is a homogeneous nonzerodivisor of degree one on $G$, we have
$\operatorname{Hilb}_{G/(x^*)}(z)=(1-z)\operatorname{Hilb}_G(z).$
	Moreover,
$	R/(x)\cong k[[Y,Z]]/(YZ,Z^2-Y^3),$
	whose $k$-basis is
$\{	1,\ y,\ y^2,\ y^3,\ z\}.$
These elements have $\mathfrak m/(x)$-adic orders $0,1,2,3,1$, respectively. Hence
$\operatorname{Hilb}_{G/(x^*)}(z)=1+2z+z^2+z^3.$
It follows that
 $\operatorname{Hilb}_G(z)
=\frac{1+2z+z^2+z^3}{1-z}=\frac{h_G(z)}{1-z}.$ 
 Recall that if $A$ is a standard graded Cohen--Macaulay ring of dimension $d$, its $a$-invariant is defined by
$a(A)=\deg h_A(z)-d$. In particular, since $G$ is one-dimensional and Cohen--Macaulay, its $a$-invariant is
$a(G)=\deg(1+2z+z^2+z^3)-1=2.$
 	
The gaps of $H$ are
$1,2,3,4,7,8,13.$
Thus $S$ is the normalization of $R$ and its conductor is
$$\mathfrak{C}:=(R:S)=u^{14}S=(u^{14},u^{15},u^{16},u^{17},u^{18})R.$$
The one-dimensional test-ideal/conductor formula gives
$\tau(R)=\mathfrak{C}$; see \cite[Theorem 3.3]{VV}. For completeness,
this formula has an elementary verification here. For every
	nonzero ideal $I\subseteq R$, the valuation criterion and the
	conductor give
$I^*=IS\cap R.$
 	For the nontrivial inclusion, choose $f\in I$ of least valuation.
 	If $w\in IS\cap R$, then $w/f\in S$; for a fixed nonzero
 	$c\in \mathfrak{C}$,
 	$cw^q=f^q c(w/f)^q\in I^{[q]}$ for every $q=p^e$.
 	Conversely, a tight-closure equation implies
 	$v(w)\geq v(I)$ by letting $q$ tend to infinity.
 	It follows that $\mathfrak{C} I^*\subseteq I$ for every $I$, so
 	$\mathfrak{C}\subseteq\tau(R)$. Taking $I=u^N R$ with $N\geq14$ gives
 	$I^*=u^N S$, and
 	$\tau(R)I^*\subseteq I$ forces $\tau(R)S\subseteq R$.
 	Hence $\tau(R)=\mathfrak{C}$.

 	The ideal $(x)$ is a reduction of $\mathfrak m$. Indeed, using the relations
 $\{yz=x^3, z^2=y^3\}$
 	we can verify that every monomial of degree four in $x,y,z$ belongs to $x\mathfrak m^3$. For the monomials not involving $x$, we have
	$$	y^4 =x\cdot xy^3,\quad
 	y^3z =x^3y^2,\quad
 	y^2z^2 =x\cdot xy^4\quad 	yz^3 =x^3z^2,\quad
 	z^4 =x\cdot xy^5 .$$ 
 	
 	Consequently,
 	$\mathfrak m^4\subseteq x\mathfrak m^3.$
 	The reverse inclusion is immediate, and hence
 	$\mathfrak m^4=x\mathfrak m^3.$
 	Therefore, $(x)$ is a reduction of $\mathfrak m$.
 	Reduction invariance and the principal multiplication formula
 	for test ideals \cite[Proposition 1.11(1)--(2)]{HY} therefore give,
 	for every integer $t\geq0$,
 	$$
 	\tau(R,\mm^t)=\tau(R,(x)^t)=x^t\tau(R)=u^{5t+14}S.
 	$$
 	Alternatively, just the containment
 	$x^t\tau(R)\subseteq\tau(R,\mm^t)$ suffices for the
 	counterexample; it follows from \cite[Lemma 2.6(1)]{WY}.
 	Since
$\mathfrak m=(u^5,u^6,u^9),$
 	every generator of $\mathfrak m^{t+3}$ is a product of $t+3$ elements, each of which has valuation at least $5$. Consequently,
 	$\mathfrak m^{t+3}
 	\subseteq u^{5(t+3)}S
 	=u^{5t+15}S.$
 	On the other hand,
 $u^{5t+14}\notin u^{5t+15}S.$
 	Therefore,
 $u^{5t+14}\notin\mathfrak m^{t+3}.$
 	Since $u^{5t+14}\in\tau(R,\mathfrak m^t)$, we conclude that
 	$
 	\tau(R,\mathfrak m^t)\not\subseteq\mathfrak m^{t+3}$
  {for every } $t\geq0$.
 	In particular, taking $t=1$, we obtain
$	u^{19}\in\tau(R,\mathfrak m)\setminus\mathfrak m^4.$ This proves the claim.
 \end{proof}

  	The proposed exponent in Question \ref{eq:q58} is too large by
 the  following:
 
 \begin{proposition}\label{prop:family}
 	For every integer $n\geq2$, let
 	$(R_n,\mm_n):=k[[u^{2n+1},u^{2n+2},u^{3n+3}]]$. Then $R_n$ is a
 	one-dimensional complete intersection, $G_n:=\gr_{\mm_n}(R_n)$
 	is Cohen--Macaulay, and
 	$$
 	a(G_n)=n,\qquad \mathfrak{C}_n=(R_n:k[[u]])=u^{2n^2+2n+2}k[[u]].
 	$$
 	For every integer $t\geq0$,
 	$\tau(R_n,\mm_n^t)=u^{et+2n^2+2n+2}k[[u]]
 	\subseteq\mm_n^{t+n}$ and
$\tau(R_n,\mm_n^t)\not\subseteq\mm_n^{t+n+1}.$

 \end{proposition}
 
 \begin{proof}Recall that the previous result is the special case $n=2$. 	Also, 
 	$	e(R)=e=2n+1$. Let $ b:=2n+2$ and $c:=3n+3.$ 
 	A presentation is given by
 $R_n\simeq k[[X,Y,Z]]/
 	(X^{n+1}-Y^{n-1}Z,\ Z^2-Y^3).$
Similar to Proposition~\ref{thm:q58-counter}, we use the
 	following basis
 	$\{1,Y,\ldots,Y^{n+1}, Z,YZ,\ldots,Y^{n-2}Z\}$  over $k[[u^e]]$.
 	Its valuations form the Ap\'ery set
 	$$
 	\{jb:0\leq j\leq n+1\}
 	\ \cup\ \{c+jb:0\leq j\leq n-2\}.
 	$$
 	These have all the distinct residues modulo $e$.
 	Their orders are respectively $0,1,\ldots,n+1$ and
 	$1,2,\ldots,n-1$, and each equals the floor of its valuation
 	divided by $e$. As before,
 	$\operatorname{ord}_{\mm_n}(u^v)=\lfloor v/e\rfloor$ for
 	$v\in\langle e,b,c\rangle$. Hence $G_n$ is Cohen--Macaulay
 	and its Hilbert numerator is
 	$$
 	1+2z+\cdots+2z^{n-1}+z^n+z^{n+1}.
 	$$
 	It follows that $a(G_n)=n$.
 	The largest Ap\'ery element is $(n+1)b$; the conductor exponent
 	is therefore $(n+1)b-e+1=2n^2+2n+2$.
 	The conductor formula and the principal reduction $(u^e)$
 	give the stated test ideal, exactly as above. Finally,
 	$ 	2n^2+2n+2=en+(n+2)$ with $ 0<n+2<e.$
 	The first element of that test ideal has order $t+n$, and
 	all of its elements have order at least $t+n$.
 \end{proof}
 
 \subsection{ 
 	$J\subseteq\tau(R,J^{d-1})$  for monomial parameter ideals}
 \label{sec:q86}
 
  \begin{question}\label{eq:q86}   (See \cite[Question 8.6]{WY}).
 	When
 	$J\subseteq\tau(R,J^{d-1}) $?\end{question}
 
 The next theorem gives a complete answer for monomial parameter
 ideals in regular parameters. The same answer holds for the
 mixed-characteristic test ideal used earlier in this note.
 
 \begin{proposition}\label{thm:q86-monomial}
 	Let $(R,\mm)$ be a complete regular  ring of dimension
 	$d\geq2$ and positive residue characteristic, and let
 	$v_1,\ldots,v_d$ be a regular system of parameters. Let
 $Q:=(v_1^{a_1},\ldots,v_d^{a_d})$ with $ a_i\geq1.$
 	In equal characteristic, assume also that $R$ is $F$-finite.
 Then
$$(\ast):\quad \quad
 		\tau(R,Q^{d-1})=
 		\left(v^\beta:\beta\in\NN^d,             \sum_{i=1}^d\frac{\beta_i+1}{a_i}>d-1\right).$$
 	In particular,
${\ Q\subseteq\tau(R,Q^{d-1})
 		\quad\Longleftrightarrow\quad
 		\sum_{i=1}^d\frac1{a_i}>d-2.\ } $
 \end{proposition}
 
 \begin{proof}
 	The Newton polyhedron of $Q$ is
 $P(Q)= \{x\in\mathbb R_{\geq0}^d:
 	\sum_i\frac{x_i}{a_i}\geq1 \}.$
 	The monomial test-ideal formula gives $(\ast)$; see
 	\cite[Theorem 4.8]{HY}, \cite[Example 2.8]{WY}, and
 	\cite[Proposition 8.21]{BMP}. For the generator $v_j^{a_j}$,
 	the strict inequality in that formula becomes
 $1+\sum_i\frac1{a_i}>d-1.$
 	It is independent of $j$, and is equivalent to the boxed
 	condition. Monomial membership in ideals generated by monomials
 	in a regular system of parameters is determined by divisibility
 	of their exponent vectors, so this checks every generator of $Q$.
 \end{proof}
 
 \begin{corollary}
If   all entries with $a_i=1$ deleted, and arrange the remaining $r$ entries
 	in increasing order, then the containment in Proposition~\ref{thm:q86-monomial}
 	holds exactly in the following cases:
 	\begin{enumerate}
 		\item $r\leq2$, with arbitrary remaining exponents;
 		\item $r=3$ and the triple is one of
 		$\{(2,2,c) ,  (2,3,3), (2,3,4), (2,3,5):c\geq2\}$
 	\end{enumerate}
Also, it never holds if $r\geq4$.
 \end{corollary}
 
 \begin{proof}
 	After removing the entries equal to $1$, the criterion is
 	$\sum_{j=1}^r1/a_j>r-2$. This is automatic for $r\leq2$.
 	For $r\geq4$, the left side is at most $r/2\leq r-2$.
 	If $r=3$ and $2\leq a\leq b\leq c$, the inequality forces
 	$a=2$. If $b=2$, it always holds; if $b=3$, it forces
 	$c=3,4,5$; and if $b\geq4$, it fails.
 \end{proof}

 \begin{corollary}
 	Adopt the notation from
 	subsection~\ref{sec:comparison}, and fix $a_1,\ldots,a_d\geq1$.
 	Set $Q_i:=(p^{a_1},X_2^{a_2},\ldots,X_d^{a_d})R_i$,
 	$Q_D:=\ulim Q_i$, and $Q^\sharp:=Q_D R^\sharp$.
 	Then the following conditions are equivalent:
 	$$
 	\begin{gathered}
 		Q_D\subseteq T_{\U}(Q^{d-1}),\qquad
 		Q^\sharp\subseteq\tau(R^\sharp,(Q^\sharp)^{d-1}),\quad
 		\sum_{i=1}^d\frac1{a_i}>d-2.
 	\end{gathered}
 	$$
 	Moreover,
 	$T_{\U}(Q^{d-1})R^\sharp=	\tau(R^\sharp,(Q^\sharp)^{d-1})$.
 \end{corollary}
 
 \begin{proof}
 In view of 	Theorem~\ref{thm:monomial}, we see the same finite monomial
 	generators in every factor and in $R^\sharp$. It remains to apply the
 	membership test in Proposition~\ref{thm:q86-monomial}.
 \end{proof}


\end{document}